\magnification1200

\font\sevenrm=tir at 7.6pt

\rm

\font\sevenbf=tib at 7.6pt

\font\tener=eurm10
\font\sevener=eurm7
\font\fiveer=eurm5
\textfont1=\tener
\scriptfont1=\sevener
\scriptscriptfont1=\fiveer
\mathcode`0="7130
\mathcode`1="7131
\mathcode`2="7132
\mathcode`3="7133
\mathcode`4="7134
\mathcode`5="7135
\mathcode`6="7136
\mathcode`7="7137
\mathcode`8="7138
\mathcode`9="7139
\font\tenfrak=eufm10

\font\tenam=msam10
\font\sevenam=msam7
\font\fiveam=msam5
\newfam\amssymfam
\textfont\amssymfam=\tenam
\scriptfont\amssymfam=\sevenam
\scriptscriptfont\amssymfam=\fiveam
\count255=\the\amssymfam
\multiply\count255 by"100
\advance\count255 by"2000
\advance\count255 by"03\mathchardef\square   =\count255
\def\qed{{\hfill$\square$}}

\font\tenbm=msbm10
\font\sevenbm=msbm7
\font\fivebm=msbm5
\newfam\bmfam
\textfont\bmfam=\tenbm
\scriptfont\bmfam=\sevenbm
\scriptscriptfont\bmfam=\fivebm
\def\Bbb#1{{\fam\bmfam #1}}

\input pst-plot

\nopagenumbers
\def\title{Real Root Conjecture fails for five and higher dimensional spheres}
\font\sevenss=cmss7
\headline{\everymath{\scriptstyle}\sevenss \ifnum\pageno>1
\ifodd\pageno\title\ \hrulefill\ \the\pageno
\else\the\pageno\ \hrulefill\ S.~R.~Gal\fi\fi} 
\centerline{\bf \title}
\smallskip
\centerline{\'Swiatos\l aw R. Gal%
\footnote{$^\star$}{Partially supported by a KBN grant 2 P03A 017 25.}}
\centerline{Wroc\l aw University}
\font\tt=cmtt10
\centerline{\tt http://www.math.uni.wroc.pl/\~{}sgal/papers/dc.ps}

\bigskip
{
\narrower\narrower\smallskip\noindent\everymath{\scriptstyle}%
\sevenbf Abstract : \sevenrm
A construction of
convex flag triangulations of five and higher dimensional spheres, whose
h-polynomials fail to have only real roots, is given.
We show that there is no such example in dimensions lower than five.

\noindent A condition weaker than realrootedness is conjectured
and some evidence is provided.\smallskip}   
\bigskip

\parindent0pt
\parskip=\smallskipamount
\footnote{}{2000 {\it Mathematics Subject Classification: }20F55.}
\footnote{}{{\it Key phrases: flag complex, h-vector, cd-index,
Charney-Davis Conjecture, Real Root Conjecture}.}

\newcount\secno\secno=0
\newcount\subsecno\subsecno=0
\newcount\thmno\thmno=0

\def\section #1\par{\vskip 0pt plus.3\vsize\penalty-250
\vskip 0pt plus-.3\vsize\bigskip\vskip\parskip
\global\advance\secno by 1
\global\subsecno=0{\tenfrak\the\secno.} {\bf #1}\par
\vskip -\medskipamount
}

\def\subsection #1\par{\vskip 0pt plus.1\vsize\penalty-250
\vskip 0pt plus-.1\vsize\medskip\vskip\parskip
\global\advance\subsecno by 1
\global\thmno=0{\tenfrak\the\secno.\the\subsecno} {\bf #1}\par
\smallskip}

\def\tag{\the\secno.\the\subsecno.\the\thmno}
\def\mktag{{\global\advance\thmno by 1}\tag}

Let the f-polynomial $f_X$
of a simplicial complex $X$ be defined by the formula
$$f_X(t)\colon=\sum_{\sigma\in X}t^{\#\sigma}.$$

There is a classical problem:
what can be said in general about the f-polynomials of
(a certain class of) simplicial complexes?

In particular, it is well known what polynomials appear as
f-polynomials of
{\parindent=.2in\parskip=0pt
\item{$\bullet$} general simplicial complexes, or
\item{$\bullet$} triangulations of spheres that are the boundary
complexes of convex polytopes
\parindent=0pt

(the reader may consult [St1] for ample discussion).
The question concerning all triangulations of spheres remains
still open. However the answer is conjecturally the same.
}

What we are interested in is the special case of the latter. Namely,
what can be said in general about the f-polynomials of
{\it flag sphere triangulations}?

The paper is organized as follows. 
In Section 1.1 we recall the definition of and basic facts about
flag complexes. In Section 1.2 we discuss the generalized homology
spheres (GHS), also called Gorenstein$^*$ complexes, a suitable
generalization of sphere triangulations.
In Section 2.1 we define a substitution in the f-polynomial of a GHS which
we call the $\gamma$-polynomial and find very useful when the complex
is flag.  In Sections 2.2 and 2.3 we discuss the relation of
$\gamma$-polynomial to the well known Charney-Davis Conjecture and
the cd-index. In Section 2.4 we explore the construction of edge
subdivision which will be used in later proofs.
In Section 3.1 we formulate the Real Root Conjecture
(a strengthening of the well known Charney-Davis Conjecture) and prove it for
generalized homology spheres of dimension at most four. In Section 3.2 we prove
a partial converse. This provides an almost complete answer to the question
posed in the previous paragraph in low dimensions. In Section 3.4
we construct a flag convex triangulation of $S^5$ that is a counterexample
to the Real Root Conjecture.

Apart from the understanding of f-polynomials of flag triangulation of
spheres the Real Root Conjecture was important since
the affirmative answer, as observed by Reiner and Welker [RW], 
would provide a partial check for the Stanley-Neggers Conjecture.

We conjecture that the coefficients of the $\gamma$-polynomial
of the flag GHS are nonnegative. We show that this conjecture is stronger
than the Charney-Davis Conjecture and weaker than the Real Root Conjecture,
thus it could remain true in all dimensions. One should think of the conjecture
on $\gamma$-polynomial as a lower bound conjecture for flag sphere
triangulations. Some evidence is provided.

The author would like to thank Tadeusz Januszkiewicz and Vic Reiner
for useful discussions and
Andrzej Derdzi\'nski and Pawe{\l} Goldstein for their extensive help with the
final version of this manuscript.

\section Preliminaries

\subsection Flag Complexes

\proclaim Definition \mktag. A simplicial complex $X$ with the vertex set $S$
is called {\rm flag} if for any $T\subset S$ such that $T$ is a {\rm clique}
(i.e., any two distinct vertices of $T$ are joined by an edge),
$T$ is a face of $X$.

Obviously a flag complex is determined by its one-skeleton.
A barycentric subdivision of any polytopial complex or, more generally,
regular CW-complex is flag.

A motivation for studying flag complexes is a theorem of Gromov,
which states that a cubical complex is locally CAT(0) if and only if
the link of any vertex is flag [Gr]. Recall that locally CAT(0) space is
in particular aspherical.

\proclaim Definition \mktag. The {\rm link} ${\rm Lk}_\sigma$
of a simplex $\sigma$ in a simplicial complex $X$ consists of all
$\tau\in X$ such that
$\sigma\cup\tau\in X$ and $\sigma\cap\tau=\emptyset$.

{\it Remark \mktag.} 
The link of any simplex in a flag complex is flag itself.

\subsection The h-polynomial

\def\GHS{{{\rm GHS}}}

\proclaim Definition \mktag. A {\rm (simplicial) generalized homology sphere
of dimension $n$ ($\GHS^n$)} is a simplicial complex such that the link of any
simplex $\sigma$ has the homology of a sphere of dimension $(n-\#\sigma)$.
We will omit the superscript if not necessary.

A simplicial GHS is also called a {\it Gorenstein$^*$} complex
(see [St1, Ch.~2, Thm.~5.2] for further reference).

\edef\ghsdef{Definition \tag}

{\it Remark \mktag.} A triangulation of a (homology) sphere
is a generalized homology sphere.

We introduce \ghsdef, since the double suspension of a GHS  is a
triangulation of a sphere [C]. Thus any GHS  may appear as the link
of an edge in some triangulation of a sphere.
More generally, assume that the dimension of
a simplicial (or cubical) complex is greater
than $2$. Then it is a cellulation of a manifold
if and only if the link of any vertex is a simply connected GHS [E,F].

\proclaim Definition \mktag. We say that $X$ is a {\rm convex sphere triangulation}
if it is the boundary complex of some convex polytope.

{\it Remark \mktag.} The link of any simplex of a convex sphere triangulation
is also a convex sphere triangulation.

\proclaim Theorem (Dehn-Sommerville relations).
If\/ $X$ is a  $\GHS^{n-1}$ then $f_X(t-1)=(-1)^nf_X(-t)$.

{\it Proof [Kl]:}
In fact, one needs a weaker assumption.
Namely, that  that $X$ is {\it Eulerian}, i.e. the Euler
characteristic of the link of any simplex $\sigma$ equals to
that of the sphere of an appropriate dimension, or in other words:
$$\sum_{\tau\supset\sigma}(-1)^{\#\tau}=(-1)^n.$$
If $\tau$ is a simplex then
$(1+s)^{\#\tau}=\sum_{\sigma\subset\tau}s^{\#\sigma}$.
Therefore
$$f_X(t-1)=\sum_\tau(t-1)^{\#\tau}=
\sum_{\sigma\subset\tau}(-t)^{\#\sigma}(-1)^{\#\tau}=
\sum_\sigma(-t)^{\#\sigma}(-1)^n=(-1)^nf_X(-t).$$
\qed 

If $X$ is a $\GHS^{n-1}$,
then a more efficient invariant is the h-polynomial defined as
$$(1+t)^nh_X\left({1\over1+t}\right)\colon=t^nf_{X}\left({1\over t}\right).$$
The Dehn-Sommerville relations written in terms of the h-polynomial
say that $h_X$ is {\it reciprocal} (i.e. $h_X(t)=t^nh_X(1/t)$).
Moreover, if $X$ is a convex
sphere triangulation then $h_X$ has a geometric interpretation in terms of a
generic height function and $h_X(t^2)$ is the Poincar\'e series of the
cohomology of the corresponding toric variety [St1].

{\it Remark \mktag.} The h-polynomial $h_X$ is usually defined for {\sl any}
simplicial complex. Although in this exposition we will use it only when
X is a GHS, the reader should be aware that some formulas (e.g. the formula
for the h-polynomial of the subdivision along an edge) are incorrect when
applied to the general complex.

\section The $\gamma$-polynomial 

\subsection The Definition

\proclaim Proposition \mktag. Assume that $h$ is a reciprocal polynomial
of degree $n$. Then there exists an unique polynomial $\gamma$ 
of degree at most $\lfloor n/2\rfloor$ with the property
$$h(t)=(1+t)^n\,\gamma\left(t\over(1+t)^2\right).\leqno(\mktag)$$
Moreover if $h$ has integral coefficients
then so does $\gamma$.

\edef\gammadef{(\tag)}

{\it Proof :} Observe that the dimension of the space of reciprocal
polynomials of degree $n$ is $\lfloor n/2\rfloor$. The polynomials
$t^i(1+t)^{n-2i}$ for $0\leq i\leq n/2$ are reciprocal and linearly
independent (being of different degree). Thus they constitute a basis
of this space. Since the leading coefficients are equal to one,
they also constitute a basis over the integers.

The coefficients of $\gamma$ are the coefficients of $h$ with respect
to this basis.\qed

{\it Remark \mktag.} If $h$ is monic then the
constant coefficient of $\gamma$ equals one.

\proclaim Definition \mktag. Let $X$ be a GHS (or at least Eulerian).
The polynomial $\gamma_X$
defined by
$$h_X(t)=(1+t)^{\deg h_X}\,\gamma_X\left(t\over(1+t)^2\right)
\leqno(\mktag)$$ will be called the $\gamma$-polynomial of $X$.

\edef\gammaxdef{(\tag)}

\proclaim Question \mktag. What is the combinatorial/geometric interpretation
of (the coefficients of) the $\gamma$-po\-ly\-no\-mial?

\proclaim Conjecture \mktag. If\/ $X$ is a flag GHS then
all the coefficients of the $\gamma$-polynomial $\gamma_X$ are nonnegative.

\edef\galconj{Conjecture \tag}

In the rest of this Section
we would like to provide some evidence for the \galconj.

\proclaim Definition \mktag. The join of two complexes $X$ and $Y$
on disjoint ground sets
is defined as the complex $\{\sigma\cup\tau\vert\sigma\in X,\ \tau\in Y\}$.

{\it Remark \mktag.} Any of $f_\bullet$, $h_\bullet$ and ${\gamma_\bullet}$
is multiplicative with respect to the joins.

\edef\multfact{Remark \tag}

\proclaim Corollary \mktag. If $X$ and $Y$ satisfy \galconj\ then so does their join
$X*Y$.

\proclaim The Generalized Lower Bound  Conjecture [St, Ch.~II, Conj.~6.2].
Let $X$ be a GHS. Then
$h_X$ is {\it unimodal}, i.e., if $h_X(t)=\sum_{i=0}^nh_it^i$
then $h_{\lfloor{n\over 2}\rfloor}\geq\dots h_2\geq h_1\geq h_0=1$.

{\it Remark \mktag.} The above conjecture is true if $X$ is a convex
triangulation of a sphere [St1].

\proclaim Corollary \mktag. If \galconj\ holds for $X$ then $h_X$ is unimodal.

{\it Proof :} Each of the polynomials $t^i(1+t)^{n-2i}$ is unimodal. A sum
of reciprocal unimodal polynomials is unimodal.\qed

\proclaim Definition \mktag. The {\rm cross polytope} $O^n$ is the $n$-fold join of
the zero-dimensional sphere.

\proclaim Lemma \mktag.
Let\/ $X$ be a flag $\GHS^{n-1}$
and $\gamma_X(t)=\sum_{i=0}^n\gamma_it^i$.
Then
{\parskip=0pt\parindent=.2in
\item{(1)} $\gamma_1\geq 0$,
\item{(2)} if\/ $\gamma_1=0$ then $X$ is a cross-polytope.
}

\edef\hlemma{Lemma \tag}

{\it Proof :}
Part (1) is equivalent to the condition that $X$ has at least
$2n$ vertices. We prove this by induction
on the dimension of $X$. Take two vertices $v$ and $w$ not joined by an edge.
The vertices adjacent to $v$ are vertices of the link of $v$,
and, by induction, there is at least $2(n-1)$ of them. Together
with $v$ and $w$ there is $2n$ of them, as desired.
This also proves (2).\qed

\subsection Charney-Davis Conjecture

\galconj\ is a strengthening of the well known
\proclaim Charney-Davis Conjecture [CD].
If $X$ is a flag $\GHS^{2n-1}$ then $(-1)^nh_X(-1)\geq0$.

Precisely, $(-1)^nh_X(-1)$
equals to the highest coefficient of $\gamma_X$. Indeed \gammaxdef\ can
be rewritten as $$\left({1\over t}\right)^nh_X(t)=
\left({(1+t)^2\over t}\right)^n\gamma_X\left({t\over(1+t)^2}\right).$$
Passing to the limit $t\to -1$ we obtain the claim.

The Charney-Davis Conjecture, obvious for $n=1$, is proven for $n=2$ [DO].
It is motivated by its being a consequence
of the Euler Characteristic Conjecture which says that the sign of the
Euler characteristic of a $2n$-dimensional closed aspherical manifold
is $(-1)^n$.

\proclaim Lemma \mktag.
The sum of the f-polynomials of the links of all simplices with
$k$ vertices of any simplicial complex $X$ is equal to $f^{(k)}_X/k!$.

{\it Proof: } For any $\tau\in{\rm Lk}_\sigma$ define $\tau^*=\tau\cup\sigma$.
We have
$$\sum_{\scriptstyle \sigma\in X\atop\scriptstyle\#\sigma=k}
f_{{\rm Lk}_\sigma}(t)=
\sum_{\scriptstyle\sigma\subset\tau^*\in X\atop\scriptstyle\mathop\#\sigma=k}
t^{\#\tau^*-k}=
\sum_{\tau^*\in X}{\#\tau^*\choose k}t^{\#\tau^*-k}={f^{(k)}_X(t)\over k!}.$$
\qed

\proclaim Corollary \mktag. If\/ $X$ is even-dimensional $GHS$ then
the highest coefficient of $\gamma_X$ is nonnegative
provided the Charney-Davis Conjecture is true for links of vertices of $X$.

\edef\evencd{Corollary \tag}

{\it Proof :} Let $\dim X=2n$. By definition
$$\sum_{i=0}^n{\gamma_i}t^i(t+1)^{2n-2i+1}=(t-1)^{2n+1}
f_X\left({1\over t-1}\right).$$

Differentiate and substitute $t=-1$. Since $f_X(-1/2)=0$ (by
the Dehn-Sommerville relations) one observes that
$${\gamma_n}=(-1)^n(1/2)^{2n-1}f_X'(-1/2)\geq0,$$
where the last inequality follows from the fact that
 $2^{2n}f_Y(-1/2)=h_Y(-1)$
for any ($2n-1$ dimensional) link $Y$ of any vertex in $X$.
\qed

\proclaim Corollary \mktag. If\/ $X$ is flag $GHS$ of dimension less than five
then the \galconj\ holds.

\edef\glowdim{Corollary \tag}

{\it Proof:} Let $(n-1)$ be the dimension of the GHS. Observe that
If $n\leq5$ then $\gamma_X$ is at most quadratic. Write
$\gamma_X(t)=1+\gamma_1t+\gamma_2t^2$.
The linear term $\gamma_1$ is nonnegative by \hlemma\ (1). Thus the claim
follows for $n\leq3$. If $n=4$ then ${\gamma_2}$ is nonnegative by the
Davis-Okun Theorem. If $n=5$ then ${\gamma_2}$ is nonnegative by the
Davis-Okun Theorem and \evencd.\qed

\subsection The cd-index

The reader may find the combinatorial background of this Section in
[St1 Ch.~III.4,St2].

Let $P$ be a {\it finite graded poset of rank $n+1$} with $\hat 0$
and $\hat 1$. Let $\rho $ be the rank function.
For any chain $C$ of the form $\hat0<x_1<\ldots<x_d<\hat1$ define 
a noncommutative monomial
$u_C=\prod_{i=1}^n u_i$ in the variables $a$ and $b$ putting
$$u_i=\cases{a&if $i\neq\rho(x_k)$ for any $k$,\cr b&if $i=\rho(x_k)$
for some $k$.\cr}$$

Finally let
$$\eqalignno{\Upsilon_P(a,b)&\colon=\sum_Cu_C\cr
\noalign{\hbox{and}}
\Psi_P(a,b)&\colon=\Upsilon_P(a-b,b).\cr}$$

One may think that $\Upsilon_P$ is a generalization
of f-polynomial, while $\Psi_P$ is that of h-polynomial
in the following sense.
Define {\it the nerve} or {\it the order complex} $N(P)$ of $P$ to be
the simplicial complex with the set of vertices $P-\{\hat0,\hat1\}$ such that
$C\subset P-\{\hat0,\hat1\}$ is in $N(P)$ if and only if $C$ is a chain. Then
$$\eqalign{f_{N(P)}(t)&=\Upsilon_P(1,t),\cr
h_{N(P)}(t)&=\Psi_P(1,t).\cr}$$
Note that $N(P)$ is always a flag complex.

\proclaim Definition \mktag. A graded poset $P$ as above is called
{\it Eulerian}
if for any $x<z$ one has $$\sum_{x\leq y\leq z}(-1)^{\rho(y)}=0.$$

\proclaim Proposition [St1,St2] (Bayer-Billera, Fine, Stanley).
If $P$ is an Eulerian poset then $\Psi_P(a,b)$
may be written as $\Phi_P(c,d)$ in $c=a+b$ and $d=ab+ba$.

The noncommutative polynomial $\Phi_P$ is called the cd-index of $P$.
It follows directly from the the definition that
$$\gamma_{N(P)}(t)=\Phi_P(1,2t).$$

{\it Remark \mktag.} E.~Babson was the first to notice that
$h_{N(P)}(-1)=\Psi_P(0,-2)$, i.e. if $n=2m$ then $(-1)^m h_{N(P)}(-1)$
equals $2^m$ times the coefficient of $d^m$ in $\Psi_P$.

The main conjecture on cd-index is the following
\proclaim Conjecture \mktag\ [St2, Conj.~2.1].
The coefficients of the cd-index
of any Gorenstein$^*$ poset are nonnegative.

\edef\stconj{Conjecture \tag}

As proven by Stanley
[St2], if $P$ is a face poset
of an S-shellable cell (e.g. a convex cell),
then \stconj\ holds for $P$.
In this case $N(P)$ is the barycentric subdivision of that cell.

\proclaim Corollary \mktag. If\/ $X$ is the barycentric subdivision of
an S-shellable cell then \galconj\ holds.

Recently Karu [Kar] have announced the proof of \stconj. As a corrolary
we obtain

\proclaim Corollary \mktag. If\/ $X$ is the barycentric subdivision of
a regular CW-sphere then \galconj\ holds.

\subsection Edge subdivision.

\def\subd{\mathop{\rm Sub}\nolimits}

\proclaim Definition \mktag. Let $X$ be a simplicial complex.
Let $\eta=\{s,t\}$ be an edge. Define $\subd_\eta(X)$ to be a simplicial
complex constructed from $X$ by bisection of all simplices containing $\eta$.
In other words let $e$ be any letter not in the vertex set $S$ of $X$.
Then $S\cup\{e\}$ is the vertex set of $\subd_\eta(X)$ and
$$\subd_\eta(X)=\{\sigma\vert\eta\not\subset\sigma\in X\}
\cup\{\sigma\cup\{e\},\sigma\cup\{s,e\},\sigma\cup\{t,e\}
\vert\sigma\in{\rm Lk}_\eta\}.$$
We say that $\subd_\eta(X)$ is a {\rm subdivision} of $X$ along $\eta$.
 
The geometric realizations of $X$ and $\subd_\eta(X)$ are
homeomorphic. In particular if $X$ is a sphere triangulation
then so is $\subd_\eta(X)$.

\proclaim Proposition \mktag. Assume that $X$ is a convex sphere triangulation.
Then $\subd_\eta(X)$ may be realized as a convex sphere triangulation.

{\it Proof: }
Take any vector $v$ starting at the midpoint of $\eta$
and pointing inside $X$.
Taking $e$ to be a sufficiently small
translation of the midpoint of $\eta$ in the direction of $-v$ we obtain that
the boundary of the convex hull of the vertices of $X$ and $e$ is a (convex)
realization of $\subd_\eta(X)$.\qed

\proclaim Proposition \mktag.
Subdividing a GHS $X$ along an edge $\eta$ affects $h$ and $\gamma$ as follows
$$\eqalign{h_{\subd_\eta(X)}(t)&=h_X(t)+t\,h_{{\rm Lk}_\eta}(t),\cr
{\gamma_{\subd_\eta(X)}}(t)&=
\gamma_X(t)+t{\gamma_{{\rm Lk}_\eta}}(t).\cr}\leqno{(\mktag)}$$

\edef\subdh{(\tag)}

{\it Proof:}
Clearly $h_{\subd_\eta(X)}-h_X$ does not depend on $X$ but only on
the link of $\eta$. Therefore it suffices to check \subdh\
in one particular case. Define $X$ to be the join of $Y$ and the $k$-gon.
Let $\eta$ be an edge of a the $k$-gon.
Then ${\rm Lk}_\eta=Y$ and $\subd_\eta(X)$ is the join of
$Y$ and the $(k+1)$-gon. Thus the claim follows from \multfact\ and
simple calculation that shows that if $\Delta_m$ is an $m$-gon, then 
$$\eqalign{f_{\Delta_m}(t)&=1+m(t+t^2),\cr
h_{\Delta_m}(t)&=1+(m-2)t+t^2,\cr
\gamma_{\Delta_m}(t)&=1+(m-4)t.\cr}\leqno{(\mktag)}$$\qed

\edef\polyh{(\tag)}

\proclaim Proposition \mktag. If\/ $X$ is flag then so is $\subd_\eta(X)$.

{\it Proof: }Let $\sigma$ be a clique in $\subd_\eta(X)$.
If $e\not\in\sigma$, then $\sigma$ is a clique in $X$.
Thus $\sigma\in X$, but then, by construction, $\sigma\in\subd_\eta(X)$.

Now assume that $e\in\sigma$. Define $\sigma_*\colon=\sigma-\{e,s,t\}$.
Each of the vertices in $\sigma_*$ being connected by an edge to $e$ has
to be connected to $s$ and $t$. Thus $\sigma_*$ is a clique in $Lk_\eta$,
and, in particular, a simplex in $X$.

On the other hand $\eta\not\subset\sigma$, since,
by definition, $\eta$ is not an edge of $\subd_\eta(X)$.
What follows $\sigma$ is of the form $\sigma_*\cup\{t,e\}$,
$\sigma_*\cup\{e\}$ or $\sigma_*\cup\{e,s\}$.
Therefore $\sigma\in\subd_\eta(X)$.\qed

\proclaim Corollary \mktag. If $X$ and $Lk_\eta$ satisfy \galconj, then
so does $\subd_\eta(X)$.

\edef\subconj{Corollary \tag}

\section Real Roots

\subsection The Conjecture 

Stating \galconj\ we were motivated by

\proclaim The Real Root Conjecture.
The zeroes of the\/ h-polynomial of a flag GHS  are all real.

{\it Remark \mktag.}
The Real Root Conjecture implies \galconj\ by the following argument. 
Assume that a reciprocal polynomial $h$ and a polynomial
$\gamma$ are related by \gammadef.
Then $h$ has only real negative roots
if and only if the same holds for $\gamma$ because
$t/(1+t)^2$ is real negative or infinite if and only if $t$ is real negative.
The ``if'' part is obvious. Conversely,
if $t/(1+t)^2=1/(\sqrt t+(\sqrt t)^{-1})^2$
is real negative then $\sqrt t+(\sqrt t)^{-1}$ is purely imaginary.
Since $\mathop{\rm Re}(z)$ and
$\mathop{\rm Re}(z^{-1})=\mathop{\rm Re}(z)/|z|^2$ have the same sign,
$\sqrt t$ has to be purely imaginary. Thus $t$ is real negative.

\edef\zeroesrem{Remark \tag}

The Real Root Conjecture was stated by T.~Januszkiewicz in a series of 
questions concerning $L^2$-cohomo\-lo\-gy of buildings (the details
may be found in [DDJO]).

Independently, V.~Reiner and V.~Welker observed that the positive answer
would prove the Neggers-Stanley Conjecture for graded naturally labeled
posets of width $2$ (details in [RW]).

{\it Remark \mktag.} As pointed out by the editor there have been recently
announced counterexamples to the Stanley-Neggers Conjecture for the
general case [Br] and for the naturally labeled one [Str].

\proclaim Theorem \mktag. The Real Root Conjecture is true
if the dimension of the GHS is less than five.

\edef\lowdim{Theorem \tag}

\proclaim Proposition \mktag.
Let\/ $h_X=\sum_{i=0}^nh_it^i$ be the h-polynomial of a flag
$\GHS^{n-1}$. Then
{\parskip=0pt\parindent=.2in
\item{(1)} among the roots of $h_X$ with the smallest modulus
there is a real negative one, and
\item{(2)}
if this root is equal to\/ $-1$, then $X$ is a cross-polytope.
}

\edef\sr{Proposition \tag}

{\it Proof :}
By the definition of the h-polynomial
$${(1+t)^n\over h_X(-t)}=\left(f_X\left({-t\over 1+t}\right)\right)^{-1}.
\leqno(\mktag)$$

\edef\pole{(\tag)}

In Section 3.3, for the sake of the self consistency of the exposition,
we will recall a proof of the following well known
\proclaim Proposition \mktag. All coefficients in the power series expansion
of \pole\ are positive.

\edef\claim{Proposition \tag}
As a corollary we obtain that there is a pole of \pole\ at
its convergence radius.
Since the modulus of any pole is not smaller than the convergence radius,
among the poles of \pole\ with the smallest modulus there is a real one.
As the poles of \pole\ have opposite values to the zeroes
of $h_X$, this proves (1).

If the root with the smallest modulus is $-1$,
then by reciprocity it has the greatest modulus, and it follows that
all roots lie on the unit circle. Thus the linear term, which is the
negative of the sum of these roots, is less than or equal to the degree.
However, \hlemma(1) gives the opposite inequality.
The equality is possible only when all roots are equal to $-1$.
Thus (2) follows from \hlemma(2).\qed

\proclaim Corollary \mktag. If\/ $X$ is not a cross-polytope, then $\gamma_X$
has a negative real root.

{\it Proof of \lowdim:}
By the \zeroesrem\ we need to check that $\gamma_X$ has only
real negative roots.
If the dimension of the flag GHS $X$ is less than five, then $\gamma_X$
is at most quadratic. As $\gamma_X$ has real coefficients
and, by previous Corollary,
at least one real root (if $\gamma_X$ is not constant), it cannot
have a single non-real root.
Since, by \glowdim, all the coefficients of $\gamma_X$
are nonnegative, it cannot have a real positive root. Thus the claim. \qed

\subsection Low dimensional geography

In this section we prove a partial converse of \lowdim.

A quadratic polynomial $1+h_1t+t^2$ has only real negative roots if and only
if $h_1\geq 2$. Then it is the h-polynomial of an $(h_1+2)$-gon.

A reciprocal polynomial $H$ of odd degree is of the form
$H(t)=(1+t)h(t)$ for some reciprocal polynomial $h$.
If $H$ has only negative roots then so does $h$.
If $h$ is the h-polynomial of some complex $X$, then $H$ is
the h-polynomial of a suspension of $X$ ($X$ joined with $S^0$). 

Therefore a reciprocal monic polynomial with integer coefficients
of degree at most three is 
the h-polynomial of a flag sphere triangulation if and only if it has
negative real roots. This is almost true if the degree is four
(or five, by the previous remark). Namely, we have

\proclaim Theorem \mktag. Let $\gamma$ be a quadratic polynomial
with constant term $1$ and integer coefficients. Assume that
$\gamma(t)-t$ has only negative real roots. Then
$h(t)=(1+t)^4\gamma(t/(1+t)^2)$ is the h-polynomial of
a flag sphere triangulation.

\edef\sfpch{Theorem \tag}

Is the condition that $\gamma(t)-t$
has only negative real zeroes essential? In other words, is there
a monic reciprocal polynomial of degree $4$ with natural coefficients
having only negative real roots such that $h$ is not the h-polynomial
of a flag sphere triangulation? The smallest example not covered by the
above Theorem is $1+9t+21t^2+9t^3+t^4$.
It is the average of the h-polynomials of
a join of a $5-$gon and an $8-$gon and a join of a $6-$gon and a $7-$gon.
However we state a

\proclaim Conjecture \mktag. Assume  that $X$ is a flag triangulation
of $S^3$ such that $\gamma_X(t)-t$ has some non-real roots.
Then $X$ is a join of two polygons.

For fixed $h_2$ there is at most one monic reciprocal polynomial
with real negative roots such that $\gamma(t)-t$ has some non-real roots.
It is not an h-polynomial of a join of polygons for $h_2=21$, $25$,
$31$, $35$, $36$, $41$, $43$, $48$, $49$, $54$, $\ldots$

\proclaim Lemma \mktag. A polynomial $h(t)=1+h_1t+h_2t^2+h_1t^3+t^4$ 
has only real negative roots if and only if 
$$\eqalign{{\bf cd}\colon=h_2-2h_1+2&\geq0,\cr
{\bf sr}\colon={h_1}^2-4(h_2-2)&\geq0,\cr
h_1&\geq 4.\cr}$$

{\it Proof: }First note that if $h(t)=(1+t)^4+{\gamma_1}t(1+t)^2
+{\gamma_2}t^2$ then by \zeroesrem\ we have to show that
$\gamma(t)=1+\gamma_1t+\gamma_2t^2$ has only real negative roots.

$\gamma$ has only real roots if an only if the discriminant
${\bf sr}={\gamma_1}\vphantom{h}^2-4{\gamma_2}$ is nonnegative.
Then the roots are negative if and only if $\gamma_2={\bf cd}$
and ${\gamma_1}=h_1-4$ are nonnegative.\qed

The plot below shows regions, in which various configurations
of roots appear. White dots mark double roots.

\vskip 0pt plus 6cm\penalty-250\vskip 0pt plus-6cm
\def\coords{
  \psline{->}(-2.2,0)(1.5,0)
  \psline{->}(0,-1.5)(0,1.5)
  \pscircle{.36}}
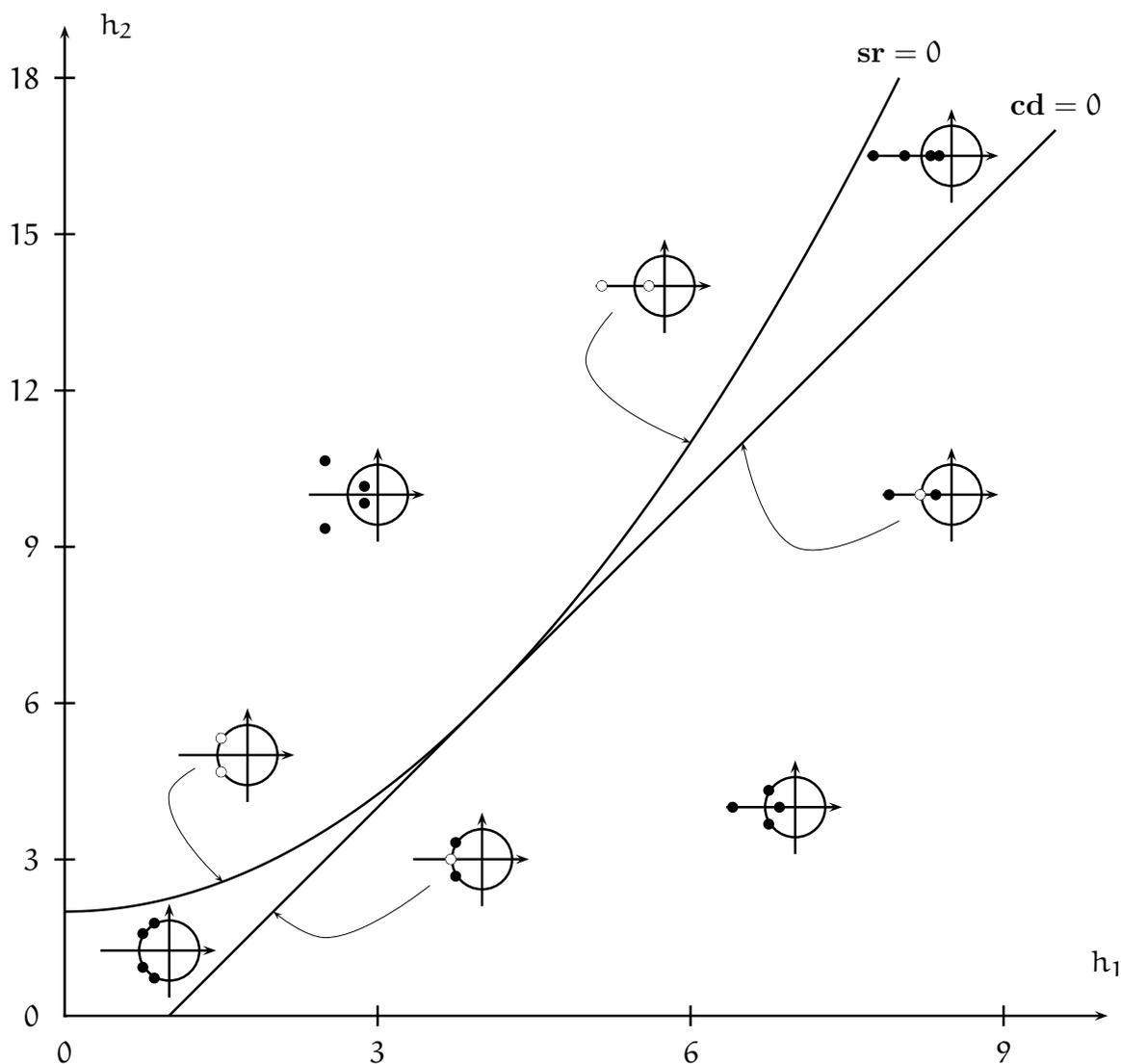
\pspicture(-1,-1)(5,12)
\psset{yunit=.6,xunit=1.2}
\psaxes[dy=3,Dy=3,dx=3,Dx=3]{->}(10,19)
\rput(10,1){$h_1$}
\rput(.5,19){$h_2$}
\psplot[plotstyle=curve]{0}{8}{x x mul 4 div 2 add}
\psline(1,0)(9.5,17)
\rput(8,18.5){${\bf sr}=0$}
\rput(9.5,17.5){${\bf cd}=0$}
\rput(7,4){\psset{xunit=.3,yunit=.6}
  \coords
  \psdots*(-.5,0)(-2,0)(-.84147,.54030)(-.84147,-.54030)}
\rput(3,10){\psset{xunit=.3,yunit=.6}
  \coords
  \psdots*(-.42735,.27015)(-.42735,-.27015)(-1.68294,1.08060)(-1.68294,-1.08060)
}
\pscurve[linewidth=.3pt]{->}(8,9.5)(7,9)(6.5,11)
\rput(8.5,10){\psset{xunit=.3,yunit=.6}
  \coords
  \psdots*(-.5,0)(-2,0)
  \psdot[dotstyle=o](-1,0)}
\pscurve[linewidth=.3pt]{->}(3.5,2.5)(2.5,1.5)(2,2)
\rput(4,3){\psset{xunit=.3,yunit=.6}
  \coords
  \psdots*(-.84147,.54030)(-.84147,-.54030)
  \psdot[dotstyle=o](-1,0)}
\pscurve[linewidth=.3pt]{->}(5.25,13.5)(5,12.5)(6,11)
\rput(5.75,14){\psset{xunit=.3,yunit=.6}
  \coords
  \psdots*[dotstyle=o](-.5,0)(-2,0)}
\pscurve[linewidth=.3pt]{->}(1.25,4.75)(1,4.25)(1.52,2.5625)
\rput(1.75,5){\psset{xunit=.3,yunit=.6}
  \coords
  \psdots*[dotstyle=o](-.84147,.54030)(-.84147,-.54030)}
\rput(1,1.25){\psset{xunit=.3,yunit=.6}
  \coords
  \psdots*(-.84147,.54030)(-.84147,-.54030)(-0.47943,-0.87758)(-0.47943,0.87758)
}
\rput(8.5,16.5){\psset{xunit=.3,yunit=.6}
  \psline{->}(-2.7,0)(1.5,0)
  \psline{->}(0,-1.5)(0,1.5)
  \pscircle{.36}
  \psdots(-.4,0)(-.666,0)(-2.5,0)(-1.5,0)}
\endpspicture

The line ${\bf cd}=0$ describes the polynomials that have zero at $-1$
(of even multiplicity, since $h$ is reciprocal). The other pair of zeroes
has to be real or lie on the unit circle. In other words the curve separates
the regions with positive and negative value of $h(-1)$ (as in
the Charney-Davis Conjecture).

The curve ${\bf sr}=0$ describes the polynomials that have
a pair of double roots. They have to be real or lie on the unit circle.
In other words the curve is the border of the region where the smallest
modulus root is real or all the roots have modulus one.

The above two curves in $(h_1,h_2)$-plane are tangent at the point $(4,6)$.
The third inequality separates two domains. In the first all roots are real,
while in the second they all lie on the unit circle.

{\it Proof of \sfpch:}
Recall that the h-polynomial is multiplicative with respect to the join and
the h-polynomial of $m$-gon equals $1+(m-2)t+t^2$.

Since $h_1=f_1-4$ is a natural number, one has
$$\left\lfloor{h_1^2\over 4}\right\rfloor=
\left\lceil{h_1\over 2}\right\rceil\left\lfloor{h_1\over 2}\right\rfloor=
\left(h_1-\left\lfloor{h_1\over 2}\right\rfloor\right)
\left\lfloor{h_1\over 2}\right\rfloor=
\max_{k\in \Bbb N} \left(h_1-k\right)k.$$

Let $\alpha_k=h_2-2-k(h_1-k)$. In particular ${\bf cd}=\alpha_2$. Define
$$C_k=\left\{(h_1,h_2)\in\Bbb Z^2\colon
\alpha_{k-1}\geq0\geq\alpha_k\right\},$$
for $k\geq 3$ and $C_2=\left\{(h_1,h_2)\in\Bbb Z^2\colon
\alpha_2=0,\ h_1\geq 4\right\}$.
Then $$\left\{(h_1,h_2)\in\Bbb Z^2\colon
6\leq2h_1-2\leq h_2\leq{h_1^2\over 4}+2\right\}=
\bigcup_{k=2}^\infty C_k.$$

$C_k$ is a cone with the vertex at the point
$(2k-1,k^2-k+2)$, which corresponds
to the h-polynomial of the join of a $(k+1)$-gon and a $(k+2)$-gon
and is generated
over $\Bbb N$ by the primitive vectors $(1,k-1)$ and $(1,k)$.

Assume that $(h_1-1,h_2-2)=(2k-1,k^2-k+2)+a(1,k-1)+b(1,k)\in C_k$.
Let $X$ be the join of a $(k+1)$-gon and
a $(k+2)-$gon subdivided along an edge, whose link is a quadrilateral.
$X$ has two disjoint edges whose links are a $k-$gon and a $(k+1)$-gon.
Subdividing the former $a$ times and the latter $b$ times we obtain a
triangulation with the desired h-polynomial.\qed

\sfpch\ should be compared to what is known about the h-polynomials of
arbitrary (not necessarily flag) sphere triangulations.
\proclaim Theorem \mktag.
A reciprocal polynomial $h(t)=\sum_{i=0}^nh_it^i$
of degree $n=4$ or $5$
is the h-poly\-nomial of a triangulation of a sphere 
if and only if $$h_1(h_1+1)/2 \geq h_2\geq h_1\geq h_0=1.$$

The ``if'' part is due to Billera and Lee [BL]. The meaning of the
inequalities is the following:
{\parindent=.2in\parskip=0pt
\item{$\bullet$} $h_1(h_1+1)/2 \geq h_2$ is equivalent to the fact
that two vertices may be joined by at most one edge,
\item{$\bullet$} $h_2\geq h_1$ is a part of the Lower Bound Theorem [Ba,Kal],
\item{$\bullet$} $h_1\geq h_0$ is equivalent to the fact that a minimum
of the number of vertices is achieved on a boundary of a simplex.
}

\subsection The Real Root Conjecture for $S^5$.

\proclaim Lemma \mktag. If\/ $1+\gamma_1t+\gamma_2t^2+\gamma_3t^3
=(1+xt)(1+yt)(1+zt)$
has only real roots then $$\gamma_2^2\geq3\gamma_3\gamma_1.\leqno(\mktag)$$

\edef\gineq{(\tag)}

{\it Proof :}
$$2\left((xy+yz+zx)^2-3xyz(x+y+z)\right)=x^2(y-z)^2+y^2(z-x)^2+z^2(x-y)^2.$$
\qed

Note that if $O$ is the cross-polytope then
${\gamma_{O}}(t)=1$.

\proclaim Theorem \mktag. Assume that $X$ is a
flag triangulation of $S^5$
that has an edge $\eta$ whose link is a cross-polytope and
$h_X(-1)<0$ (cubical coefficient of $\gamma_X$ is positive).
Then for a sufficiently large natural number $m$, the
$m$-fold subdivision $\subd^m_\eta(X)$
of\/ $X$ along $\eta$ contradicts the Real Root Conjecture, i.e.
$h_{\subd^m_\eta(X)}$ has a non-real root. 

\edef\ceth{Theorem \tag}

{\it Proof :}
By \subdh, subdivision along $\eta$ increases the linear term of $\gamma_X$
without changing other coefficients. But \gineq\ for $\gamma_X$ is false when
$\gamma_3=-h_X(-1)$ is positive and the linear coefficient $\gamma_1$
is sufficiently large.\qed

What is left to do is to exhibit a complex satisfying the hypothesis of \ceth.

Take the join $X_1$ of two pentagons.
Subdivide an edge whose link is a quadrilateral to obtain $X_2$.
Finally let $X$ be the join
of $X_2$ and another pentagon. $X_1$ has a vertex (the new one) whose link is
a cross polytope. Thus, $X$ has an edge whose link is a cross polytope
(join of the above with any edge of the pentagon).

By \subdh\ and \polyh\ we calculate:
$$\gamma_X(t)=(1+t)((1+t)^2+t\cdot1)=1+4t+4t^2+t^3.$$
Therefore $X$ satisfies the hypothesis of \ceth\ and, twice subdivided, $X$
becomes a counterexample to the Real Root Conjecture.
In fact, it suffices to take
a single subdivision, but this requires an extra check.
This (smallest konown) counterexample to the Real Root Conjecture
has f-vector $f_{\subd_\eta X}(t)=1+17\,t+109\,t^2+
345\,t^3+575\,t^4+483\,t^5+161\,t^6$.

The above example still satisfies \galconj\
and, in particular, the Charney-Davis Conjecture (cf. \subconj).

\proclaim Corollary \mktag. Taking the join with any flag sphere triangulations
one finds that there are counterexamples to the Real Root Conjecture
if the dimension of the sphere is greater or equal to $5$.

Other counterexamples to the Real Root Conjecture are presented in
the forthcoming paper [Ga].

\subsection The smallest root.

We briefly present the classical proof of \claim\ that uses Coxeter groups.
Another homological proof can be given by showing
that Stanley-Reisner face ring $R(X)$ of
a flag complex $X$ is Koszul, thus the coefficient of \pole\  
at $t^j$ equals $\dim Tor^{R(X)}_j(k,k)$ (see e.g. [RW, Prop. 4.13]).

With any flag complex $X$ with vertex set $S$ one associates the
{\it right angled Coxeter group} $W$ with the following  presentation
$$W=\langle S\vert s^2=1\hbox{\rm\ for all }s\in S,
st=ts\hbox{\rm\ for all }\{s,t\}\in X\rangle.$$
The definition of $W$ uses only the one-skeleton of $X$,
but the following observation links the whole $X$ to $W$.
Let the subgroup of $W$ generated by $T\subset S$
be denoted by $W_T$.  Then $W_T$ is finite if
and only if $T\in X$. In this case $W_T=\left(\Bbb Z/2\right)^{\#T}$.

\proclaim Definition \mktag. Define a formal series
$W(t)=\sum_{w\in W}t^{\ell(w)}$,
where $\ell$ denotes the length function with respect to the
generating set $S$. We call $W(\cdot)$
the growth series of $W$. 

\proclaim Proposition \mktag\ [Se]. $W(t)$ represents a rational function.
Moreover, if $W$ is infinite,
then $${1\over W(t)}=\sum_{T\subset S}{(-1)^{\#T}\over W_T(t^{-1})},$$
where $T$ runs over subsets of $S$ such that $W_T$ is finite.

\edef\serre{Proposition \tag}

\proclaim Corollary \mktag.
Let $W$ be the Coxeter group associated to a flag complex $X$.
Then
$$f_{X}\left({-t\over1+t}\right)
={1\over W(t)}. \leqno(\mktag)$$
\global\edef\fw{(\tag)}

{\it Proof: }
If $W_T$ is finite, then $W_T(t)=(1+t)^{\#T}$, thus \serre\ reduces to \fw.\qed

This finishes the proof of \claim.

Note that proving \sr(1)\ we did not use the assumption that $X$ is GHS.
The second part of \sr\ allows the following generalization
(one may consult [DDJO, Prop.~3.10] for ample discussion):

\proclaim Proposition \mktag. If $X$ is any flag complex and the radius of
convergence of $W_X(\cdot)$ equals one (i.e. $h_X$ has no zeroes
in the interior of the unit disk), then $X$ is a join of
a cross-polytope and a simplex (i.e. a multiple suspension of a simplex).

{\it Proof :} Let $T$ be any subset of $S$.
The coefficients of $W_T(\cdot)$ are dominated by those of $W(\cdot)$.
This is straightforward either by noticing that the length function
on $W_T$ is the restriction of the length function on $W$ or
by interpreting the coefficients as dimensions of Tor modules.

What follows, the convergence radius of the $W_T(\cdot)$
is greater or equal to that of $W(\cdot)$.

Let $x$ and $y$ be two vertices not joined by an edge
(otherwise $X$ is a simplex and we are done).
We want to show that $X$ is a suspension of its subcomplex
spanned by $S-\{x,y\}$. To do this we need to show that if
$z$ is any vertex different from $x$ and $y$
then $z$ is joined with both $x$ and $y$. 

If not, then straightforward computation, using \serre,
shows that the convergence radius
of $W_{\{x,y,z\}}(\cdot)$ is strictly smaller than $1$.\qed

\bigskip{\bf References}\par\smallskip

\parindent=.4in
\item{[Ba]} D.~Barnette,
{\it A proof of the lower bound conjecture for convex polytopes},
Pacific J. Math. {\bf 46} (1973), pp.~349--354.
\item{[Br]} P.~Br\"and\'en, {\it Counterexamples to the Neggers-Stanley
Conjecture}, arXiv:math.CO/0408312,
\item{[BL]} L.~Billera and C.~Lee,
{\it Sufficiency of McMullen's conditions for $f$-vectors of simplicial
polytopes},
Bull. Amer. Math. Soc. {\bf 2} (1980), pp.~181--185.
\item{[C]} J.~W.~Cannon, {\it Shrinking cell-like decompositions of manifolds.
Codimension three.}, Ann.\ Math. {\bf 110} (1979), pp.~83--112
\item{[CD]} R.~Charney and M.~Davis, {\it The Euler characteristic
of a nonpositively curved, piecewise Euclidean manifold},
Pacific J.~Math. {\bf 171} (1995), pp.~117--137,
\item{[DDJO]} M.~Davis, J.~Dymara, T.~Januszkiewicz and B.~Okun,
{\it Decompositions of Hecke - von Neumann modules and the
$L^2$-cohomology of buildings}, arXiv:math.GT/0402377,
\item{[DO]} M.~Davis and B.~Okun, {\it Vanishing theorems and conjectures
for the $L^2-$homology of right-angled Coxeter Groups}, Geom.~Topol.
{\bf 5} (2001), pp.~7--74,
\item{[E]} R.~D.~Edwards, {\it The topology of manifolds and cell-like maps},
Proc. ICM Helsinki, 1978, pp.~111-127,
\item{[F]}M.~H.~Freedman, {\it The topology of four-dimensional manifolds}, 
J.~Diff.~Geom.  {\bf 17}  (1982), no. 3, pp.~357--453.
\item{[Ga]} S.~R.~Gal, {\it On Normal Subgroups of Coxeter Groups
Generated by Standard Parabolic Subgroups}, 2004 preprint,
\item{[Gr]} M.~Gromov, {\it Hyperbolic Groups}, in Essays in Group Theory,
S.~G.~Gersten ed., Springer Verlag, MSRI Publ. {\bf 8} (1987), pp.~75--263,
\item{[Kal]} G. Kalai, {\it  Rigidity and the lower bound theorem.~I.},
Invent. Math. {\bf 88} (1987), no.~1, pp.~125--151
\item{[Kar]} K.~Karu, {\it The cd-index of fans and lattices},
arXiv:math.AG/0410513,
\item{[Kl]} V.~Klee, {\it A combinatorial proof of Poincar\'e's duality
theorem}, Can.~J.~Math. {\bf 16} (1964), pp.~517--531,
\item{[RW]} V.~Reiner and V.~Welker, {\it On the Charney-Davis
and Neggers-Stanley Conjectures}, 2002 preprint,
\item{[Se]} J.~P.~Serre, {\it Cohomologie des groupes discrets}, in
{\it Prospects in Mathematics}, pp.~77--169,
Annal of Math.~Studies No. {\bf 70}, Princeton 1971,
\item{[St1]} R.~Stanley, {\it Combinatorics and commutative algebra},
Birkh\"auser 1996,
\item{[St2]} R.~Stanley, {\it Flag f-vectors and the cd-index},
Math.~Z. {\bf 216} (1994), pp.~483--499,
\item{[Ste]} J.~Stembridge, {\it Counterexamples to the Poset Conjecture
of Neggers, Stanley, and Stembridge}, 2004 preprint.
\bigskip

\rightline{\it Wroc\l aw, 31 December 2003/1 January 2004}
\bye